\documentclass[a4paper,11pt,centertags,psamsfonts]{amsart}
\usepackage{amsmath,latexsym}
\usepackage{amssymb}

\usepackage{times}

\newcommand{\ii}{\mathrm{i}}

\DeclareSymbolFont{SY}{U}{psy}{m}{n}
\DeclareMathSymbol{\emptyset}{\mathord}{SY}{'306}
\newcommand{\bA}{\mathbf A}
\newcommand{\bB}{\mathbf B}
\newcommand{\bbR}{\mathbb R}

\newcommand{\bV}{\mathbf V}

\renewcommand{\epsilon}{\varepsilon}
\newcommand{\spec}{\mathrm{spec}}
\newcommand{\C}{\mathbb{C}}

\newcommand{\R}{\mathbb{R}}
\newcommand{\EE}{\mathsf{E}}

\newcommand{\cB}{{\mathcal B}}

\newcommand{\cG}{\mathcal{G}}

\newcommand{\cM}{{\mathcal M}}

\renewcommand{\Im}{{\ensuremath{\mathrm{Im}}}}

\newcommand{\dist}{{\ensuremath{\mathrm{dist}}}}

\newcommand{\tr}{\mathrm{tr}}

\newcommand{\fH}{\mathfrak{H}}

\newcommand{\fM}{\mathfrak{M}}
\newcommand{\fN}{\mathfrak{N}}

\DeclareMathOperator{\Ran}{\mathrm{Ran}}

\newcommand{\Dom}{\mathrm{Dom}}

\newtheorem{theorem}{Theorem}[section]
\newtheorem{lemma}[theorem]{Lemma}
\newtheorem{proposition}[theorem]{Proposition}

\newtheorem{remark}[theorem]{Remark}
\newtheorem{hypothesis}[theorem]{Hypothesis}
\newtheorem{definition}[theorem]{Definition}
{\bf}{\rm}

\newtheorem{introtheorem}{Theorem}{\bf}{\it}
{\bf}{\it}
{\bf}{\it}
{\bf}{\it}

\begin{document}

\title[Non-Closed Invariant Subspaces]{The Singularly Continuous Spectrum
and Non-Closed Invariant Subspaces}

\author[V.~Kostrykin]{Vadim Kostrykin}
\address{Fraunhofer-Institut f\"{u}r Lasertechnik \\
Steinbachstra{\ss}e 15 \\ Aachen, D-52074, Germany}
\email{kostrykin@ilt.fraunhofer.de,
kostrykin@t-online.de
\newline{\rm URL}
\tt http://home.t-online.de/home/kostrykin}

\author[K.~A.~Makarov]{Konstantin A.~Makarov}
\address{Department of Mathematics\\ University of Missouri\\
Co\-lum\-bia, MO 65211, USA} \email{makarov@math.missouri.edu \newline{\rm URL}
\tt http://www.math.missouri.edu/people/kmakarov.html}

\dedicatory{Dedicated to Israel Gohberg on the occasion of his 75-th
birthday}

\begin{abstract}
Let $\mathbf{A}$ be a bounded self-adjoint operator on a separable Hilbert
space $\mathfrak{H}$ and $\mathfrak{H}_0\subset\mathfrak{H}$ a closed
invariant subspace of $\mathbf{A}$. Assuming that $\mathfrak{H}_0$ is of
codimension $1$, we study the variation of the invariant subspace
$\mathfrak{H}_0$ under bounded self-adjoint perturbations $\mathbf{V}$ of
$\mathbf{A}$ that are off-diagonal with respect to the decomposition
$\mathfrak{H}= \mathfrak{H}_0\oplus\mathfrak{H}_1$. In particular, we prove
the existence of a one-parameter family of dense non-closed invariant
subspaces of the operator $\bA+\bV$ provided that this operator has a
nonempty singularly continuous spectrum. We show that such subspaces are
related to non-closable densely defined solutions of the operator Riccati
equation associated with generalized eigenfunctions corresponding to the
singularly continuous spectrum of $\bB$.
\end{abstract}

\subjclass{Primary 47A55, 47A15; Secondary 47B15}

\keywords{Invariant subspaces, operator Riccati equation, singular
spectrum.}

\maketitle

\section{Introduction}\label{Intro}

In the present article we address the problem of a perturbation of
invariant subspaces of self-adjoint operators on a separable Hilbert space
$\fH$ and related questions on the existence of solutions to the operator
Riccati equation.

Given a self-adjoint operator $\bA$ and a closed invariant subspace
$\fH_0\subset\fH$ of $\bA$ we set $A_i=\bA|_{\fH_i}$, $i=0,1,$ with
$\fH_1=\fH\ominus\fH_0$. Assuming that the perturbation $\bV$ is
off-diagonal with respect to the orthogonal decomposition $\fH=\fH_0 \oplus
\fH_1$ consider the self-adjoint operator
\begin{equation*}
\bB = \bA +
\bV=\begin{pmatrix} A_0 & V \\ V^\ast & A_1
\end{pmatrix}\quad\text{with}\quad\bV=\begin{pmatrix} 0 & V \\ V^\ast & 0
\end{pmatrix},
\end{equation*}
where $V$ is a linear operator from $\fH_1$ to $\fH_0$. It is well known
(see, e.g., \cite{Kostrykin:Makarov:Motovilov:2}) that the Riccati equation
\begin{equation}\label{intro:Ricca}
A_1 X - X A_0 -X V X +
V^\ast =0
\end{equation}
has a closed (possibly unbounded) solution $X:\fH_0 \rightarrow\fH_1$ if
and only if its graph
\begin{equation}\label{intro:Graph}
\cG(\fH_0,X) := \left\{x\in \fH\,|\,x=x_0\oplus X x_0,\,
x_0\in\Dom(X)\subset\fH_0 \right\}
\end{equation}
is an invariant closed subspace for the operator
$\bB$.

Sufficient  conditions guaranteeing the existence of a solution to
equation \eqref{intro:Ricca} require in general  the assumption
that the spectra of the operators $A_0$ and $A_1$ are
separated,
\begin{equation}\label{separated}
d:=\dist(\spec(A_0), \spec(A_1)) > 0,
\end{equation}
and hence $\fH_0$ and
$\fH_1$ are necessarily spectral invariant subspaces of the operator $\bA$. In
particular (see \cite{Kostrykin:Makarov:Motovilov:5}), if
\begin{equation}\label{V:small}
\|V\| < c_\pi d\quad \text{with}\quad c_\pi=\frac{3 \pi - \sqrt{\pi^2 +
32}}{\pi^2 - 4}=0.503288\ldots,
\end{equation}
then the Riccati equation \eqref{intro:Ricca} has a bounded solution $X$ satisfying the bound
\begin{equation*}
\frac{\|X\|}{\sqrt{1+\|X\|^2}}\leq\frac{\pi}{2}
\frac{\|V\|}{d-\delta_V}<1
\end{equation*}
with
\begin{equation*}
\delta_V=\|V\| \tan\bigg (\frac{1}{2} \arctan\frac{2\|V\|}{d}\bigg ).
\end{equation*}
It is plausible to conjecture that  condition \eqref{V:small} can
be relaxed by the weaker requirement $\|V\|<\sqrt{3} d/2$ (see
\cite{Kostrykin:Makarov:Motovilov:5} for details). However, no
proof of that is available as yet.

In general, without additional assumptions, neither condition
\eqref{separated} nor a smallness assumption like \eqref{V:small} on the
magnitude of the perturbation $V$ can be dropped. However, if the spectra of
$A_0$ and $A_1$ are subordinated in the sense that
\begin{equation*}
\sup\spec(A_0)\leq \inf\spec(A_1),
\end{equation*}
then for any $V$ with arbitrary large norm the Riccati equation
\eqref{intro:Ricca} has a contractive solution
\cite{Kostrykin:Makarov:Motovilov:3} (see also
\cite{Adamyan:Langer:Tretter:2000a}). Note that in this case the invariant
subspaces $\fH_0$ and $\fH_1$ are not necessarily supposed to be spectral
invariant subspaces of $\bA$.

In the present work we prove  new existence results for the Riccati
equation under the assumption that the subspace $\fH_1$ is
\emph{one-dimensional}. In particular, these results imply the existence of
a one-parameter family of non-closed invariant subspaces of the
self-adjoint operator $\bB$, provided that $\bB$ has nonempty singularly
continuous spectrum.

The main result of our paper is presented by the following theorem.

\begin{introtheorem}\label{intro:sc:pp}
Assume that $\dim\fH_1=1$ and suppose that $\fH_0$ is a cyclic subspace for
the operator $A_0$ generated by the one-dimensional subspace $\Ran V$. Let
$S_{\mathrm{pp}}$ denote the set of all eigenvalues of the operator $\bB$.

Then there exists  a minimal support $S_{\mathrm{s}}$ of the
singular part of the spectral measure of the operator $\bB$ such
that:

(i) For any $\lambda\in S_{\mathrm{sc}}=S_{\mathrm{s}}\setminus
S_{\mathrm{pp}}$ the subspace $\Psi(\lambda)=\cG(\fH_0,X_\lambda)\subset\fH$
is a dense non-closed graph subspace
with $X_\lambda:\fH_0\rightarrow\fH_1$ a non-closed densely
defined operator solving the Riccati equation \eqref{intro:Ricca}
in the sense of Definition \ref{DefRicc} below.

(ii) For any $\lambda\in S_{\mathrm{pp}}\subset S_{\mathrm{s}}$
the subspace $\Psi(\lambda)=\cG(\fH_0,X_\lambda)\subset\fH$ is a closed
graph subspace  of codimension $1$ with
$X_\lambda:\fH_0\rightarrow\fH_1$ a bounded operator solving the
Riccati equation \eqref{intro:Ricca}. Moreover, the operator
$X_\lambda$ is an isolated point (in the operator norm topology)
of the set of all bounded solutions to the
Riccati equation.

The mapping $\Psi$ from $S_{\mathrm{s}}$ to the set $\cM(\bB)$ of
all (not necessarily closed) subspaces of $\fH$ invariant with
respect to the operator $\bB$ is injective.
\end{introtheorem}

The article is organized as follows. In Section \ref{sec:nonclosed} we
establish a link  between non-closable densely defined solutions to the
Riccati equation \eqref{intro:Ricca} and the associated non-closed invariant
subspaces of the operator $\bB$. In Section \ref{sec:2} accommodating the
Simon-Wolff theory \cite{Simon:Wolff} to rank two off-diagonal perturbations we
perform the spectral analysis of this operator under the assumption that
$\dim\fH_1=1$. The main result of this section is Theorem \ref{thm:x}.
Theorem \ref{intro:sc:pp} will be proven in Section \ref{sec:Riccati}.

Throughout the whole work the Hilbert space $\fH$ will assumed to be
separable. The notation $\cB(\fM,\fN)$ is used for the set of
bounded linear operators from the Hilbert space $\fM$ to the Hilbert space
$\fN$. We will write $\cB(\fN)$ instead of $\cB(\fN,\fN)$.

\section{Non-Closed Graph Subspaces}\label{sec:nonclosed}
  Let $\fH_0$ be a closed subspace
of a Hilbert space $\fH$ and $X$ a densely defined (possibly
unbounded and not necessarily closed) operator from $\fH_0$ to
$\fH_1=\fH_0^\perp:=\fH\ominus\fH_0$ with domain $\Dom(X)$. A
linear subspace
\begin{equation*}
\cG(\fH_0,X) := \left\{x\in \fH\,|\,x=x_0\oplus X x_0,\,
x_0\in\Dom(X)\subset\fH_0 \right\}
\end{equation*}
is called the graph subspace of $\fH$ associated with the pair $(\fH_0,X)$ or, in
short, the graph of $X$.

Recalling general facts on densely defined closable operators
(see, e.g., \cite{Kato:book}) we mention the
following: If $X:\; \fH_0\rightarrow\fH_1$ is a
densely defined non-closable operator, then $\cG(\fH_0,X)$ is a
non-closed subspace of $\fH$. Its closure is not a graph subspace,
i.e., there is no closed operator $Y$ such that
\begin{equation*}
\overline{\cG(\fH_0,X)}=\cG(\fH_0,Y).
\end{equation*}

\begin{proposition}\label{cor:4.2}
Let $X:\; \fH_0\rightarrow\fH_1$ be a densely defined non-closable
operator. Then the closed subspace $\overline{\cG(\fH_0,X)}$ contains an
element orthogonal to $\fH_0$.
\end{proposition}

\begin{proof}
First, for $X:\; \fH_0\rightarrow\fH_1$ being a densely defined
non-closable operator we prove the following alternative: either the closed
subspace $\overline{\cG(\fH_0,X)}$ contains an element orthogonal to
$\fH_0$ or the subspace $\fH_0$ contains an element orthogonal to
$\overline{\cG(\fH_0,X)}$. Indeed, assume on the contrary that neither the
closed subspace $\overline{\cG(\fH_0,X)}$ contains an element orthogonal to
$\fH_0$ nor the subspace $\fH_0$ contains an element orthogonal to
$\overline{\cG(\fH_0,X)}$. Then by Theorem 3.2 in
\cite{Kostrykin:Makarov:Motovilov:2} there is a closed densely defined
operator $Y:\; \fH_0\rightarrow\fH_1$ such that
$\overline{\cG(\fH_0,X)}=\cG(\fH_0,Y)$, which is a contradiction.

Now assume that the subspace $\fH_0$ contains an element $x_0$ orthogonal
to $\overline{\cG(\fH_0,X)}$. Obviously, this element is orthogonal to
$\cG(\fH_0,X)$, that is, $\langle x_0\oplus 0, x_0\oplus X x_0\rangle=0$,
and hence $x_0=0$. Then, by the alternative proven above the subspace
$\overline{\cG(\fH_0,X)}$ contains an element orthogonal to $\fH_0$,
completing the proof.
\end{proof}

For notational setup assume the following  hypothesis.

\begin{hypothesis}\label{hyp0}
 Let $\bB$ be a self-adjoint operator
represented with respect to the decomposition $\fH=\fH_0 \oplus
\fH_1$ as a $2\times2$ operator block matrix
\begin{equation}\label{1}
 \bB=\begin{pmatrix}
  A_0 & V \\
  V^\ast & A_1
\end{pmatrix},
\end{equation}
where $A_i\in\cB(\fH_i)$, $i=0,1$, are bounded self-adjoint operators in
$\fH_i$ while $V\in\cB(\fH_1,\fH_0)$ is a bounded operator from $\fH_1$ to
$\fH_0$. More explicitly, $\bB =\bA+\bV$, where $\bA$ is the bounded
diagonal self-adjoint operator,
\begin{equation}\label{2}
 \bA = \begin{pmatrix}
  A_0 & 0 \\
  0 & A_1
\end{pmatrix},
\end{equation}
and the operator $\bV=\bV^\ast$ is an off-diagonal bounded operator
\begin{equation}\label{3}
 \bV=\begin{pmatrix}
  0 & V \\
  V^\ast & 0
\end{pmatrix}.
\end{equation}
\end{hypothesis}

\begin{definition}\label{DefRicc}
A densely defined (possibly unbounded and not necessarily closable)
operator $X$ from $\fH_0$ to $\fH_1$ with domain $\Dom (X)$ is called a
\emph{strong solution} to the Riccati equation
\begin{equation}\label{Riccati:2} A_1 X - X A_0 - XVX + V^\ast =0
\end{equation}
if
\begin{equation*}
\left.\Ran(A_0+VX)\right|_{\Dom(X)}\subset\Dom(X)
\end{equation*}
and
\begin{equation*}
A_1 X x - X (A_0 + VX) x + V^\ast x =0\quad\text{for any}\quad x\in\Dom(X).
\end{equation*} \end{definition}

\begin{theorem}\label{2x2}
Assume Hypothesis \ref{hyp0}. A densely
defined (possibly unbounded and not necessarily closed) operator $X$ from
$\fH_0$ to $\fH_1$ with domain $\Dom(X)$ is a strong solution to the
Riccati equation \eqref{Riccati:2} if and only if the graph subspace
$\cG(\fH_0, X)$ is invariant for the operator $\bB$.
\end{theorem}

\begin{proof}
First, assume that $\cG(\fH_0, X)$ is invariant for $\bB$. Then
\begin{equation*}
\bB(x\oplus Xx)=(A_0 x + V X x)\oplus(A_1 X x + V^\ast x)\in\cG(\fH_0, X)
\end{equation*}
for any $x\in \Dom (X)$. In particular, $A_0 x + V X x\in \Dom(X)$ and
\begin{equation*}
A_1 X x + V^\ast x = X (A_0 x + V X x) \text{ for all } x\in \Dom (X),
\end{equation*}
which proves that $X$ is a strong solution to the Riccati
equation \eqref{Riccati}.

To prove the converse statement assume that $X$ is a strong solution to the
Riccati equation \eqref{Riccati:2}, that is,
\begin{equation*}
A_0 x + V X x
\in \Dom (X)
\end{equation*}
and
\begin{equation*}
A_1 X x + V^\ast x = X
(A_0 x+ V X x),\qquad x\in \Dom (X),
\end{equation*}
which proves that the
graph subspace $\cG(\fH_0, X)$ is $\bB$-invariant.
\end{proof}

\begin{remark}
By Lemma 4.3 in \cite{Kostrykin:Makarov:Motovilov:2} a closed densely
defined operator $X:\; \fH_0\rightarrow\fH_1$ is a strong solution to the
Riccati equation \eqref{Riccati:2} if and only if it is a weak solution to
\eqref{Riccati:2}.
\end{remark}

\section{The Singular Spectrum of the Operator $\bB$}\label{sec:2}

Assume the following hypothesis.

\begin{hypothesis}\label{hyp} Assume Hypothesis \ref{hyp0}.
Assume in addition that the
Hilbert space $\fH_1$ is one-dimensional,
\begin{equation*}
\fH_1=\C,
\end{equation*}
and the Hilbert space $\fH_0$ is the cyclic subspace generated by $\Ran V$.
\end{hypothesis}

Note that under  Hypothesis \ref{hyp} the Hilbert space $\fH_0$
can be realized as a space of square integrable functions with
respect to a Borel probability measure $m$ with
compact support,
\begin{equation*}
\fH_0 = L^2(\R;m)
\end{equation*}
such that the bounded operator $A_0$ acts on $L^2(\R,m)$ as the
multiplication operator
\begin{equation*}
(A_0 x_0)(\lambda)=\lambda x_0(\lambda),\qquad x_0\in L^2(\R,m),
\end{equation*}
$A_1$ is the multiplication by a real number $a_1$
and, finally, the linear bounded map
\begin{equation*}
V^\ast\;:\; \fH_0\rightarrow\fH_1
\end{equation*}
is given by
\begin{equation*}
V^\ast x_0 = \langle v, x_0\rangle_{\fH_0},\qquad x_0\in\fH_0
\end{equation*}
for some $v\in\fH_0$.

\begin{lemma}\label{lem:simple}
Assume Hypothesis \ref{hyp}. Then the element $0\oplus
1\in\fH=\fH_0\oplus\fH_1$ is cyclic for the operator $\bB$ given by
\eqref{1} -- \eqref{3} and, hence, $\bB$ has a simple spectrum.
\end{lemma}

\begin{proof}
By hypothesis (in the above notations) the element
$v\in\fH_0$ is cyclic for the operator $A_0$. Therefore, the
cyclic subspace with respect to the operator $\bB$ generated by
the elements $v\oplus 0\in\fH$ and $0\oplus 1\in\fH$ is the whole
$\fH$. Without loss of generality we may assume that $a_1=0$.
Observing that $\bB(0\oplus 1) = v\oplus 0$ proves the claim.
\end{proof}

\begin{theorem}\label{equiv:spec}
Assume Hypothesis \ref{hyp}.
Then the Herglotz function
\begin{equation}\label{Stern}
\phi(z) = \frac{1 + (a_1 - z) \langle v, (A_0 - z)^{-1} v\rangle_{\fH_0}}
{(a_1 - z) - \langle v, (A_0 - z)^{-1} v\rangle_{\fH_0}}
\end{equation}
admits the representation
\begin{equation*}
\phi (z)
=\int \frac{d\omega(\lambda)}{\lambda-z},
\end{equation*}
where $\omega$ is a probability measure on $\R$ with compact
support. Moreover, the operator $\bB$ is unitarily equivalent to
the multiplication operator by the independent variable on
$L^2(\bbR, \omega)$.
\end{theorem}

\begin{proof}
Introduce the Borel measure $\Omega$ with values in the set of non-negative
operators on $\fH_1\oplus\fH_1$ by
 \begin{equation*}
\Omega(\delta)=\begin{pmatrix} V & 0 \\ 0 & 1 \end{pmatrix}^\ast
\EE_{\bB}(\delta)
\begin{pmatrix} V & 0 \\ 0 & 1 \end{pmatrix},
\end{equation*}
where $\begin{pmatrix} V & 0 \\ 0 & 1
\end{pmatrix}$ is the linear map from $\fH_1\oplus\fH_1$ to
$\fH_0\oplus\fH_1$  and let
\begin{equation*}
\omega(\delta)=\tr\,\Omega(\delta),\qquad \delta\subset\R\text{ a Borel
set}.
\end{equation*}
Clearly, the measure $\omega$ vanishes on all Borel sets $\delta$ such that $\EE_{\bB}(\delta)=0$. In fact,
these measures have the same families of Borel sets, on which they vanish.
Indeed, assuming
$\omega(\delta)=0$ yields
\begin{equation*}
\langle v\oplus 0, \EE_{\bB}(\delta)\, v\oplus 0\rangle_{\fH} + \langle
0\oplus 1, \EE_{\bB}(\delta)\, 0\oplus 1\rangle_{\fH} = 0
\end{equation*}
and, hence, in particular,
\begin{equation}\label{nuli}
\langle 0\oplus 1, \EE_{\bB}(\delta)\, 0\oplus 1\rangle_{\fH} = 0,
\end{equation}
which  implies  $\EE_{\bB}(\delta)=0$.

Introducing  the $\cB(\fH_1\oplus\fH_1)$-valued Herglotz function
\begin{equation}\label{M:def}
M(z) = \begin{pmatrix} V & 0 \\ 0 & 1
\end{pmatrix}^\ast(\bB-z)^{-1} \begin{pmatrix} V & 0 \\ 0 & 1
\end{pmatrix}
\end{equation}
one concludes that the Herglotz function $M(z)$ admits the representation
\begin{equation*}
M(z)=\int_\R \frac{d\Omega(\lambda)}{\lambda-z},
\end{equation*}
and hence
\begin{equation*}
\tr M(z)=\int_\R \frac{d\omega(\lambda)}{\lambda-z}.
\end{equation*}

Straightforward computations show that the operator-valued function
\eqref{M:def} with respect to the orthogonal decomposition $\fH=\fH_0\oplus
\fH_1$ can be represented as  the $2\times 2$ matrix
\begin{equation*} {M}(z) = \begin{pmatrix} {M}_{00}(z) & {M}_{01}(z)
\\ {M}_{10}(z) & {M}_{11}(z)\end{pmatrix}
\end{equation*}
with the  entries given by \begin{equation*} \begin{split}
{M}_{00}(z) &= (a_1-z)\langle v, (A_0-z)^{-1}v\rangle [a_1-z-\langle v, (A_0-z)^{-1}v\rangle]^{-1},\\
{M}_{11}(z) &= [a_1-z-\langle v, (A_0-z)^{-1}v\rangle]^{-1},\\
{M}_{01}(z) &= -(a_1-z)^{-1}{M}_{00}(z),\\
{M}_{10}(z) &= - (a_1-z)^{-1}{M}_{00}(z).
\end{split}
\end{equation*}
Taking the trace of $M(z)$ yields representation \eqref{Stern}.

Since by Lemma \ref{lem:simple} the element $0\oplus 1 $ is cyclic
and the measure $\omega$ and the spectral measure
$\EE_{\bB}$ have the same families of Borel sets, on which they vanish, one concludes (see,
e.g., \cite{Birman:Solomyak}) that the operator $\bB$ is unitarily
equivalent to the multiplication operator by the independent variable on
$L^2(\R, \omega)$, completing the proof.
\end{proof}

Recall that a measurable not necessarily closed set $S\subset\R$ is a
support of a measure $\nu$ if $\nu(\R\setminus S)=0$. A support $S$ is said
to be minimal if any measurable subset $S^\prime\subset S$ with
$\nu(S^\prime)=0$ has Lebesgue measure zero.

\begin{theorem}\label{thm:x}
The sets
\begin{equation}\label{s}
S_{\mathrm{s}} := \left\{\lambda\in\R\Big|\ a_1-\lambda
=\int\frac{|v(\mu)|^2 d m(\mu)}{\mu-\lambda-\ii 0} \right\}
\end{equation}
 and
\begin{equation}\label{sc}
S_{\mathrm{sc}} := \left\{\lambda\in\R\Big|\ a_1-\lambda =
\int\frac{|v(\mu)|^2 d m(\mu)}{\mu-\lambda-\ii 0},\quad \int\frac{|v(\mu)|^2
d m(\mu)}{|\mu-\lambda|^2}=\infty\right\}
\end{equation}
are minimal supports of the singular part $\omega_{\mathrm{s}}$ and the
singularly continuous part $\omega_{\mathrm{sc}}$ of the measure $\omega$,
respectively. The set
\begin{equation}\label{pp}
S_{\mathrm{pp}} := \left\{\lambda\in\R\Big|\ a_1-\lambda =
\int\frac{|v(\mu)|^2 d m(\mu)}{\mu-\lambda},\quad \int\frac{|v(\mu)|^2 d
m(\mu)}{|\mu-\lambda|^2}<\infty\right\}
\end{equation}
coincides with the set of all atoms of the measure $\omega$.
\end{theorem}

\begin{proof}
The fact that \eqref{s} is a minimal support of $\omega_{\mathrm{s}}$
follows from  Lemma 3.5 in \cite{Gilbert}, where one sets
$m_a^+(z)=(a_1-z)$ and
\begin{equation*}
m_b^+(z)=\langle v, (A_0-z)^{-1}
v\rangle_{\fH_0}= \int\frac{|v(\mu)|^2 d m(\mu)}{\mu-z}, \quad \Im \, z\ne 0.
\end{equation*}

It is not hard to see (cf., e.g., Example 1 in \cite{Albeverio})
 that the set $S_{\mathrm{pp}}$ coincides with the set of all
eigenvalues of the operator $\bB$. Hence, by Theorem
\ref{equiv:spec} one proves that $S_{\mathrm{pp}}$ coincides with the set
of all atoms of the measure $\omega$. Therefore,
to prove that \eqref{sc} is a minimal support of $\omega_{\mathrm{sc}}$ it
suffices to check the inclusion
\begin{equation}\label{inclus}
S_{\mathrm{pp}}\subset S_{\mathrm{s}}.
\end{equation}

Assume that  $\lambda\in S_{\mathrm{pp}}$, that is,
\begin{equation}\label{pass}
 a_1-\lambda=\int\frac{|v(\mu)|^2 d
m(\mu)}{\mu-\lambda}
\end{equation}
and
\begin{equation*}\label{fini}
\int\frac{|v(\mu)|^2 d m(\mu)}{|\mu-\lambda|^2} < \infty.
\end{equation*}
Since
\begin{equation*}
\int\frac{|v(\mu)|^2 d
m(\mu)}{|\mu-\lambda|}\le \bigg (\int\frac{|v(\mu)|^2 d
m(\mu)}{|\mu-\lambda|^2}\bigg )^{1/2}\|v\|_{L^2(\R;m)},
\end{equation*}
the dominated convergence theorem yields
\begin{equation*}
\int\frac{|v(\mu)|^2 d m(\mu)}{\mu-\lambda-\ii 0}\equiv\lim_{\varepsilon
\to +0}\int\frac{|v(\mu)|^2 d m(\mu)}{\mu-\lambda-\ii \varepsilon}
=\int\frac{|v(\mu)|^2 d m(\mu)}{\mu-\lambda},
\end{equation*}
which together with \eqref{pass} proves inclusion \eqref{inclus}. The proof is
complete.
\end{proof}

\begin{remark}
By Lemma 5 in \cite{Gilbert:Pearson} from Theorem \ref{equiv:spec} it
follows that there exist minimal supports of the absolutely continuous part
$\omega_{\mathrm{ac}}$, the singular part
$\omega_{\mathrm{s}}$, and the singularly continuous part
$\omega_{\mathrm{sc}}$ of the measure $\omega$ such that their closures
coincide with the absolute continuous part
$\spec_{\mathrm{ac}}(\bB)$, the singular part $\spec_{\mathrm{s}}(\bB)$,
and the singularly continuous part $\spec_{\mathrm{sc}}(\bB)$ of the
spectrum, respectively.
\end{remark}

\section{Riccati Equation}\label{sec:Riccati}

Given $\lambda\in\R$, introduce the  operator (linear functional)
\begin{equation*}
X_\lambda: L^2(\R;m)\to \fH_1=\C
\end{equation*}
on
\begin{equation*}
\Dom(X_\lambda)=\left\{\varphi\in L^2(\R;m)\, \Big|
\lim_{\epsilon\rightarrow +0}
\int\frac{\overline{v(\mu)}\varphi(\mu)}{\mu-\lambda-\ii \epsilon}
dm(\mu)\,\,\,\text{ exists  finitely}\right\}
\end{equation*}
by
\begin{equation}\label{Xzeta:def}
X_{\lambda}\varphi = \lim_{\epsilon\rightarrow +0}
\int\frac{\overline{v(\mu)}\varphi(\mu)}{\mu-\lambda-\ii \epsilon}
dm(\mu),\qquad \varphi\in \Dom(X_\lambda).
\end{equation}

\begin{lemma}\label{beschraenkt}
If $\lambda\in S_{\mathrm{s}}$, then the operator $X_{\lambda}$ is densely
defined.
\end{lemma}

\begin{proof}
Since the element
$v\in L^2(\R;m)$ is generating for the operator $A_0$, the set
\begin{equation*}
D=\{\varphi\,\,\, \vert \,\,\varphi(\mu)=v(\mu)\psi(\mu),
\,\,\psi\text{ is continuously differentiable on  } \R\}
\end{equation*}
is dense in $L^2(\R;m)$. For $\varphi \in D$  and $\varepsilon >0$ one obtains
\begin{align}\label{oz}
\int\frac{\overline{v(\mu)}\varphi(\mu)}{\mu-\lambda-\ii \epsilon} dm(\mu)
&=\psi(\lambda)\int\frac{|v(\mu)|^2}{\mu-\lambda-\ii \epsilon}
dm(\mu)\\\label{int} &+\int\frac{|v(\mu)|^2(\psi(\mu)-
\psi(\lambda))}{\mu-\lambda-\ii \epsilon} dm(\mu).
\end{align}
Since $\lambda\in S_{\mathrm{s}}$, by Theorem \ref{thm:x} the limit
\begin{equation*}
\lim_{\varepsilon \to +0}\int\frac{|v(\mu)|^2 d m(\mu)}{\mu-\lambda-\ii
\varepsilon}=\int\frac{|v(\mu)|^2 d m(\mu)}{\mu-\lambda-\ii 0}
\end{equation*}
exists  finitely. The integral  \eqref{int} also has a limit as
$\varepsilon \to +0$ since $\psi$ is a continuously differentiable which
proves that the left hand side of \eqref{oz} has a finite limit as
$\varepsilon \to +0$. Therefore, $D\subset \Dom(X_\lambda)$, that is,
$X_{\lambda}$ is densely defined.
\end{proof}

\begin{remark}\label{RR}
Note that by the Riesz representation theorem $X_\lambda$ is
bounded whenever the condition
\begin{equation}\label{iii}
\int\frac{|v(\mu)|^2}{|\lambda-\mu|^2}dm(\mu)<\infty
\end{equation}
holds true. The converse is also true: If $X_\lambda$ is bounded, then
\eqref{iii} holds. Indeed, by the uniform boundedness principle from
definition \eqref{Xzeta:def} it follows that
\begin{equation*}
\sup_{\epsilon\in(0,1]}\int\frac{|v(\mu)|^2}{(\mu-\lambda)^2+\epsilon^2}
dm(\mu) <\infty,
\end{equation*}
proving \eqref{iii} by the monotone convergence theorem.
\end{remark}

\begin{theorem}\label{Ric:main}
Let $\lambda\in S_{\mathrm{s}}$. Then the operator $X_{\lambda}$ is a
strong solution to the Riccati equation
\begin{equation}\label{Riccati}
A_1 X - X A_0 - XVX + V^\ast = 0.
\end{equation}
Moreover, if $\lambda\in S_{\mathrm{pp}}$, the solution $X_\lambda$ is
bounded and if $\lambda\in S_{\mathrm{sc}}= S_{\mathrm{s}}\setminus
S_{\mathrm{pp}}$, the operator $X_\lambda$ is non-closable.
\end{theorem}

\begin{proof}
Note that $A_0\, \Dom(X_\lambda)\subset \Dom(X_\lambda)$. If $\lambda\in
S_{\mathrm{s}}$, then by Theorem \ref{thm:x}
\begin{equation*}
a_1-\lambda = \int\frac{|v(\mu)|^2 d m(\mu)}{\mu-\lambda-\ii 0}.
\end{equation*}
In particular, $v\in \Dom(X_\lambda)$ and
\begin{align*}
X_{\lambda}VX_\lambda\varphi&= \int\frac{|v(\mu)|^2}{\mu-\lambda-\ii 0}
dm(\mu) \cdot X_\lambda\varphi \\&= (a_1-\lambda)
\int\frac{\overline{v(\mu)}\varphi(\mu)}{\mu-\lambda-\ii 0} dm(\mu),\quad
\varphi\in \Dom(X_\lambda).
\end{align*}
Therefore, for an arbitrary $\varphi\in \Dom(X_\lambda)$ one gets
\begin{align*}
&A_1 X_\lambda \varphi - X_\lambda A_0\varphi - X_\lambda V
X_\lambda\varphi \\
&= \int \frac{\overline{v(\mu)}\varphi(\mu)(a_1-\mu)}{\mu-\lambda-\ii 0}
  dm(\mu)-(a_1-\lambda)
\int\frac{\overline{v(\mu)}\varphi(\mu)}{\mu-\lambda-\ii 0}
dm(\mu)\\
& = \int\frac{\overline{v(\mu)}\varphi(\mu)(\lambda-\mu)}{\mu-\lambda-\ii
0} dm(\mu) = -\int \overline{v(\mu)}\varphi(\mu) dm(\mu) =- V^\ast\varphi,
\end{align*}
which proves that the operator $X_{\lambda}$ is a strong solution to the
Riccati equation \eqref{Riccati}.

If  $\lambda\in S_{\mathrm{pp}}$, then \eqref{iii} holds, in which
case $X_\lambda$ is bounded. If $\lambda\in
S_{\mathrm{sc}}=S_{\mathrm{s}}\setminus S_{\mathrm{pp}}$,
then $X_\lambda$ is an
unbounded densely defined operator (functional) (cf. Remark
\ref{RR}). Since every closed finite-rank operator is bounded
\cite{Kato:book}, it follows that for $\lambda\in S_{\mathrm{sc}}$
the unbounded solution $X_\lambda$ is non-closable.
\end{proof}

\begin{proof}[Proof of Theorem \ref{intro:sc:pp}]
Introduce the mapping
\begin{equation}\label{psi}
\Psi(\lambda) = \cG(\fH_0,X_{\lambda}), \quad \lambda\in S_{\mathrm{s}},
\end{equation}
where $X_{\lambda}$ is the strong solution to the Riccati equation referred
to in Theorem \ref{Ric:main}. By Theorem \ref{2x2} the subspace
$\Psi(\lambda)$, $\lambda\in S_{\mathrm{s}}$
is invariant with respect to $\bB$. To prove the injectivity
of the mapping $\Psi$, assume that $\Psi(\lambda_1)=\Psi(\lambda_2)$ for
some $\lambda_1,\lambda_2\in S_{\mathrm{s}}$. Due to \eqref{psi},
$X_{\lambda_1}=X_{\lambda_2}$ which by \eqref{Xzeta:def} implies
$\lambda_1=\lambda_2$.

(i). Let $\lambda\in S_{\mathrm{sc}}$. By Theorem \ref{Ric:main} the
functional $X_\lambda$ is non-closable. Since $X_\lambda$ is densely
defined, the closure $\overline{\cG(\fH_0, X_\lambda)}$ of the subspace
$\cG(\fH_0, X_\lambda)$ contains the subspace $\fH_0$. By Proposition
\ref{cor:4.2}, the subspace $\overline{\cG(\fH_0, X_\lambda)}$ contains an element orthogonal
to $\fH_0$. Since $\fH_0\subset\fH$ is of codimension $1$, one concludes
that $\overline{\cG(\fH_0, X_\lambda)}=\fH_0\oplus\fH_1=\fH$.

(ii). Let $\lambda\in S_{\mathrm{pp}}$. By Theorem 5.3 in
\cite{Kostrykin:Makarov:Motovilov:2} the solution $X_{\lambda}$ is an
isolated point (in the operator norm topology) of the set of all bounded
solutions to the Riccati equation \eqref{Riccati} if and only if the
subspace $\cG(\fH_0,X_{\lambda})$ is spectral, that is, there is a Borel set
$\Delta\subset\R$ such that
\begin{equation*}
\cG(\fH_0,X_{\lambda}) = \Ran \EE_{\bB}(\Delta).
\end{equation*}

Observe that the one-dimensional graph subspace
$\cG(\fH_1,-X_{\lambda}^\ast)$ is invariant with respect to the operator
$\bB$. This subspace is spectral since by Lemma \ref{lem:simple} $\lambda$
is a simple eigenvalue of the operator $\bB$. Thus,
$\cG(\fH_0,X_{\lambda})=\cG(\fH_1,-X_{\lambda}^\ast)^\perp$ is also a
spectral subspace of the operator $\bB$.
\end{proof}

\subsection*{Acknowledgments}
The authors are grateful to C.~van der Mee for useful
suggestions. K.~A.~Makarov is indebted to the Graduiertenkolleg
``Hierarchie und Symmetrie in mathematischen Modellen'' for its kind
hospitality during his stay at the RWTH Aachen in the Summer of 2003.


\end{document}